\documentstyle[12pt]{article}
\expandafter\chardef\csname pre amssym.def
at\endcsname=\the\catcode`\@
\catcode`\@=11

\def\undefine#1{\let#1\undefined}
\def\newsymbol#1#2#3#4#5{\let\next@\relax
 \ifnum#2=\@ne\let\next@\msafam@\else
 \ifnum#2=\tw@\let\next@\msbfam@\fi\fi
 \mathchardef#1="#3\next@#4#5}
\def\mathhexbox@#1#2#3{\relax
 \ifmmode\mathpalette{}{\m@th\mathchar"#1#2#3}%
 \else\leavevmode\hbox{$\m@th\mathchar"#1#2#3$}\fi}
\def\hexnumber@#1{\ifcase#1 0\or 1\or 2\or 3\or 4\or 5\or 6\or 7\or 8\or
 9\or A\or B\or C\or D\or E\or F\fi}

\font\tenmsa=msam10 \@magscale1 \font\sevenmsa=msam7 \@magscale1
\font\fivemsa=msam5 \@magscale1
\newfam\msafam
\textfont\msafam=\tenmsa \scriptfont\msafam=\sevenmsa
\scriptscriptfont\msafam=\fivemsa
\edef\msafam@{\hexnumber@\msafam} \mathchardef\dabar@"0\msafam@39
\def\dashrightarrow{\mathrel{\dabar@\dabar@\mathchar"0\msafam@4B}}
\def\dashleftarrow{\mathrel{\mathchar"0\msafam@4C\dabar@\dabar@}}

\def\ulcorner{\delimiter"4\msafam@70\msafam@70 }
\def\urcorner{\delimiter"5\msafam@71\msafam@71 }
\def\llcorner{\delimiter"4\msafam@78\msafam@78 }
\def\lrcorner{\delimiter"5\msafam@79\msafam@79 }
\def\yen{{\mathhexbox@\msafam@55 }}
\def\checkmark{{\mathhexbox@\msafam@58 }}
\def\circledR{{\mathhexbox@\msafam@72 }}
\def\maltese{{\mathhexbox@\msafam@7A }}

\font\tenmsb=msbm10 \@magscale1 \font\sevenmsb=msbm7 \@magscale1
\font\fivemsb=msbm5 \@magscale1
\newfam\msbfam
\textfont\msbfam=\tenmsb \scriptfont\msbfam=\sevenmsb
\scriptscriptfont\msbfam=\fivemsb
\edef\msbfam@{\hexnumber@\msbfam}
\def\Bbb#1{\fam\msbfam\relax#1}
\font\teneufm=eufm10 \@magscale1 \font\seveneufm=eufm7 \@magscale1
\font\fiveeufm=eufm5 \@magscale1
\newfam\eufmfam
\textfont\eufmfam=\teneufm \scriptfont\eufmfam=\seveneufm
\scriptscriptfont\eufmfam=\fiveeufm

\catcode`\@=\csname pre amssym.def at\endcsname

\newcommand{\kla}{\left ( }
\newcommand{\mer}{\right ) }

\newcommand{\for}{\begin{eqnarray*}}
\newcommand{\mel}{\end{eqnarray*}}

\newcommand{\mitt}{\left | { \atop } \right.}
\newcommand{\kl}{\pl \le \pl}
\newcommand{\gl}{\pl \ge \pl}

\newcommand{\lel}{\pl = \pl}

\newcommand{\nz}{{\rm  I\! N}}

\newcommand{\rz}{{\Bbb R}}

\newcommand{\ten}{\otimes}

\newcommand{\p}{\hspace{.05cm}}
\newcommand{\pl}{\hspace{.1cm}}
\newcommand{\pll}{\hspace{.3cm}}

\newcommand{\hz}{\vspace{0.5cm}}

\newcommand{\qed}{\hspace*{\fill}$\Box$\hz\pagebreak[1]}

\newcommand{\Om}{\Omega}
\newcommand{\om}{\omega}

\newcommand{\al}{\alpha}
\newcommand{\si}{\sigma}

\newcommand{\la}{\lambda}
\newcommand{\eps}{\varepsilon}

\newcommand{\U}{{\cal U}}

\newcommand{\noo}{\left \|}
\newcommand{\rrm}{\right \|}

\newcommand{\intt}{\int\limits}
\newcommand{\summ}{\sum\limits}

\newcommand{\C}{{\mathcal C}}

\newtheorem{lemma}{Lemma}[section]
\newtheorem{prop}[lemma]{Proposition}
\newtheorem{theorem}[lemma]{Theorem}
\newtheorem{cor}[lemma]{Corollary}

\newtheorem{rem}[lemma]{Remark}

\newcommand{\re}{\begin{rem}\rm}
\newcommand{\mar}{\end{rem}}
\newtheorem{exam}[lemma]{Example}

\begin{document}
\title{ \bf The optimal order for the p-th moment
of sums of independent random variables with respect to symmetric
norms and related combinatorial estimates }
\author{{\large \bf  Marius Junge}\\
{\small University of Illinois at Urbana-Champaign} \\[-0.2cm]
{\small 1409 West Green Street}\\[-0.2cm]   {\small Urbana, Illinois 61801,USA} \\[-0.2cm]
{\small  e-mail: junge@math.uiuc.edu} }
\date{}
\maketitle

\begin{abstract} For $n$ independent random variables $f_1,..,f_n$ and a symmetric norm $\noo \pl \rrm_X$ on $\rz^n$,
we show that for $1\le p<\infty$
 \for
 \lefteqn{ \frac{1}{2+4\sqrt{2}} \kla
 (n\intt_0^{\frac 1n} h^*(s)^p ds )^{\frac{1}{p}} +
 \noo \summ_{i=1}^n  (n\intt_{\frac{i-1}{n}}^{\frac
 in} h^*(s) ds ) \pl e_i \rrm_X \mer  }\\
 & & \kl \kla \intt \noo \summ_{i=1}^n f_i e_i\rrm_X^p d\mu
 \mer^{\frac{1}{p}} \\
 & & \kl \frac{c_0 p}{1+\ln p} \pl \kla
  (n\intt_0^{\frac 1n} h^*(s)^p ds )^{\frac{1}{p}} +
 \noo \summ_{i=1}^n  (n\intt_{\frac{i-1}{n}}^{\frac
 in} h^*(s) ds ) \pl e_i \rrm_X \mer
 \pl .
 \mel
Here
\[ h(t,\om) \lel \summ_{i=1}^n 1_{[\frac{i-1}{n},\frac{i}{n})}(t) f_i(\om) \]
is the disjoint sum of the $f_i$'s and $h^*$ is the
non-increasing rearrangement. Similar results (where $L_p$ is
replaced by a more general rearrangement invariant function
space) were obtained first by Litvak, Gordon, Sch\"{u}tt and
Werner for Orlicz spaces $X$ and independently by S.
Montgomery-Smith for general $X$ but without an explicit analysis
of the order of growth for the constant in the upper estimate. The
order $\frac{p}{1+\ln p}$ is optimal and obtained from
combinatorial estimates for doubly stochastic matrices. The
result extends to Lorentz-norms $l_{f,q}$ on $\rz^n$ under mild
assumptions on $f$. We give applications to the theory of
noncommutative $L_p$ spaces.
\end{abstract}

\newpage

{\Large \bf Introduction and Notation} \hz

The interaction between Banach space theory and probabilistic
methods have a long tradition and Rosenthal's inequality
\cite{Ro}, extended by Burkholder \cite{Bu} to general
martingales, is an example for an inequality motivated by Banach
space theory with a significant impact in probability. Since then
there has been a big progress in calculating the expectation of
the norm of independent variables in particular by Johnson,
Schechtman, Zinn and  Johnson, Schechtman \cite{JS},  Kwapien,
Szulga \cite{KSz}, Hitczenko \cite{Hit},
Hitczenko/Montgomery-Smith \cite{HMS} and many others. Motivated
by embedding problems for non-commutative $L_p$ spaces, we
calculate the $p$-norm of the sum of independent random variables
with respect to symmetric norms\footnote{The first version of
this article dates back to 1998}. This extends recent results of
Gordon, Litvak, Sch\"{u}tt and Werner. Very recently, similar
results have been obtained by Montgomery-Smith \cite{MS} (but by
different techniques and without an analysis of the order of
constant  involved. The starting point of our approach an article
of Geiss \cite{Gei} using certain $K$-functional from
interpolation theory. \hz

More precisely, we calculate  the $p$-th moment of the sum of $n$
independent random variables $f_1,..,f_n$ with respect to  a
symmetric norm $\noo \pl \rrm_X$ on $\rz^n$, i.e. a norm
satisfying
\[ \noo \summ_{i=1}^n \eps_i x_{\pi(i)} e_i\rrm_X \lel \noo \summ_{i=1}^n x_i e_i\rrm_X \pl \]
for all coefficients $\al_i$, changes of signs $\eps_i\pm1$ and
permutations $\pi:\{1,..,n\}\to \{1,..,n\}$. Our results follow a
general philosophy: {\it Independent variables behave like
disjoint variables}. Although our results hold for all such
symmetric norms, it turns out that only few  classical norms are
really relevant for this investigation. Indeed, certain
$K$-functionals $|\pl \_k$  between the $\ell_1$ and the
$\ell_\infty$ norm (see section 1) and the weak-$\ell_1$ norm
$\ell_{1,\infty}$ (see section 2 and section 3). \hz

Let us recall the notation of the non-increasing rearrangement.
Given $n$ independent random variables $f_1,..,f_n$, we consider
the disjoint sum
\[  h(t,\om) \lel \summ_{i=1}^n 1_{[\frac{i-1}{n},\frac{i}{n})}(t) f_i(\om)\pl . \]
Then non-increasing rearrangement $h^*$ is defined by
\[ h^*(s) \lel \inf\{ t \p|\p Prob(|h|>t)\le s\} \pl .\]
A starting point of  our approach  is the following theorem of S.
Geiss \cite[proof of theorem 3.4]{Gei} based on previous work of
Johnson and Schechtman \cite{JS}.

\begin{theorem}[S. Geiss] \label{GEI} Let $f_1,..,f_n$ be independent ranodm variables, then
\[ 2^{-\frac{1}{p}} \pl (n\intt_0^{\frac{1}{n}} h^*(s)^p ds )^{\frac{1}{p}}
\kl  \kla \intt_\Om \sup_{i}|f_i|^p d\mu \mer^{\frac{1}{p}} \kl
2^{1-\frac{1}{p}} \pl(n\intt_0^{\frac{1}{n}} h^*(s)^p ds
)^{\frac{1}{p}}   \pl .\]
\end{theorem}

On the other hand, we are motivated by some combinatorial
estimates of Kwapien and Sch\"utt \cite{KSI,KSII}. In contrast to
our situation they consider
\[ \intt_{\pi \mbox{ \scriptsize Permutation}}  \noo \summ_{i=1}^n \al_{i\pi(i)} \rrm_X d\pi \]
for the $\ell_p$-norms or more generally Orlicz spaces and
arbitrary matrices $\al_{ij}$. This is also the starting point of
\cite{GSW}. In general results for the permutation group or
independent coefficients are very similar. We refer to
\cite{MS-S} for further information on averages with respect to
the group of permutations. Our main result is the following.

\begin{theorem} \label{main}  Let $1\le p <\infty$ and $\noo \pl \rrm_X$ be a
$1$-symmetric, $1$-unconditional  norm on $\rz^n$ with
normalized  unit vectors $e_i$. Then for all independent random
variables $f_1,..,f_n$ on a probability space $(\Om,\mu)$
 \for
 \lefteqn{  \frac{1}{2+4\sqrt{2}}
   \kla (n\intt_0^{\frac{1}{n}}
 h^*(s)^p ds)^\frac{1}{p} + \noo \summ_{i=1}^n
 (n\intt_{\frac{i-1}{n}}^{\frac{i}{n}}
 h^*(s) \p ds ) e_i\rrm_X  \mer
  \kl  \kla \intt_\Om \noo
 \summ_{i=1}^n f_i e_i\rrm_X^p d\mu \mer^{\frac1p} }\\
 & &\kl
  \frac{c_0p}{1+\ln p} \pl
  \kla (n\intt_0^{\frac{1}{n}}
 h^*(s)^p ds)^\frac{1}{p} + \noo \summ_{i=1}^n
 (n\intt_{\frac{i-1}{n}}^{\frac{i}{n}}
 h^*(s) \p ds ) e_i\rrm_X  \mer \pl . \hspace{3cm} {\atop }
 \mel
The order of the constant   $\frac{p}{1+\ln p}$ is optimal.
\end{theorem}

For the lower estimate, it  suffices  to consider the supremums
norm and the sequence of  $|\pl |_k$-norms (see section 1). In
this paper we give a selfcontained proof of the second  inequality
(called upper bound) only using Rosenthal's inequality and new
combinatorial tools. It turns out that an upper bound the norm in
Lorentz space $\ell_{1,\infty}$ implies the general case (see
section 2 for this.)  Indeed, the key estimate is of
combinatorial nature and we believe it is of independent interest.

\begin{theorem} \label{comb} There exists a constant $c_0$ with the
following property. Let $(\mu_{ij})_{ij} $ be a doubly
stochastic  matrix, i.e.for all $i,j$
\[ \summ_k \mu_{ik} \lel 1\lel \summ_k \mu_{kj} \pl .\]
Then
 \[ \kla \summ_{j_1,..,j_n=1}^n \kla \sup_r \frac{1}{r} {\rm card}\{ i \p|\p j_i\le r\}
 \mer^p \prod_{k=1}^n \mu_{kj_k}  \mer^{\frac{1}{p}}
 \kl c_0 \frac{p}{1+\ln p} \pl .\]
The order of growth is optimal.
\end{theorem}

As an application of our result, we obtain a result for the
Schatten class $S_X$ associated to a symmetric sequence space
 \[ S_X \lel \left\{a\in B(\ell_2) \mitt  \noo
 a\rrm_{S_X} \lel \noo \summ_k s_k(a)e_k\rrm_X
 \pl<\pl \infty  \right\} \pl .\]
Here $s_k(a)=\la_k(\sqrt{a^*a})$ are the singular values of $a$.

\begin{theorem} \label{m4} Let $1\le p<\infty$ and $X$ be a symmetric
sequence space. Then there is a von Neumann algebra $N$ such
that  $X$ embeds into some $L_p(N)$ iff there is a (possibly
different) von Neumann algebra $N$ such that $S_X$ embeds into
$L_p(N)$.
\end{theorem}

We will use standard notation from probability. $|A|$ denotes the
cardinality of a set. Moreover, $c_0$ is used (as above) for an
absolute constant varying in each occurrence. We use $a\sim_c b$
if $\frac{1}{\sqrt{c}}\kl \frac{a}{b} \le \sqrt{c}$. The space of
sequences converging to $0$ is denoted by $c_o$. The paper is
organized as follows. We prove the lower estimate in section 1.
In the second part, we show how the upper estimate can be deduced
from Rosenthal's inequality and the combinatorial estimate.
Section 4 is devoted to the elementary proof of the combinatorial
estimate. Theorem \ref{m4} and further applications are contained
in section 5. \hz

{\bf Acknowledgment:} I  want to thank Stefan Geiss for helpful
discussions leading to the conjecture of theorem \ref{main} and
Y. Gordon for bringing to his attention the articles \cite{GSW}
and \cite{MS}.

\section{The lower estimate}
In this section, we will prove the lower estimate in   Theorem
\ref{main}.  The arguments are simialr as in \cite{GSW} and
\cite{MS}, but we give a proof with a `concrete' estimate. We
assume that $\noo \pl \rrm_X$ is a $1$-unconditional,
$1$-symmetric and normalized norm. In our investigation, the
following norms are of particular interest
 \for
 \noo x\rrm_{p\pll {\atop} } &=& \kla \summ_{i=1}^n |x_i|^p \mer^{\frac{1}{p}} \pl ,\\
 \noo x\rrm_{p,\infty} &=& \sup_{k} k^{\frac{1}{p}} \pl x^*_k \pl,\\
 |x|_{k\pll \p {\atop}} &=& \summ_{j=1}^k x^*_j \pl .
 \mel
Here $(x^*_j)_{j=1}^n$ given by
 \[ x_j^*\pl:=\pl \inf\left\{ t \mitt |\{i\p|\p |x_i|>t\}|\pl <\pl j\right\} \]
is the non-increasing rearrangement of $x$. Indeed, $|\pl|_k$ is
equivalent to the $K$-functional between $\ell_1$ and
$\ell_\infty$ at the value $k$. We will need the following lemma
of \cite{KSI} which is the analogue of Geiss' theorem \ref{GEI}.

\begin{theorem}[Kwapien Sch\"utt] \label{KS} Let $(\al_{ij})_{ij=1,..,n}$ be an $n\times n$ matrix, $\al^*_1,..,\al^*_n$
be its non-decreasing rearrangement and $d\pi$ the normalized
counting measure on the group $\Pi_n$ of permutations of
$\{1,..,n\}$, then
\[ \frac{1}{2} \pl \frac{1}{n} \summ_{j=1}^n \al^*_j \kl
\intt_{\Pi_n} \sup_{i=1,..,n} |\al_{\pi(i)i}|\pl d\pi \kl
\frac{1}{n} \summ_{j=1}^n \al^*_j  \pl .\]
\end{theorem}

\begin{cor} \label{kk}  Let $1\le k\le n$ and $x=(x_1,..,x_n)\in \rz^n$, then
\[ \frac{1}{4} \pl \frac{1}{k} \summ_{j=1}^k x^*_j \kl
\intt_{\Pi_n} \sup_{i\le \frac{n}{k} } |x_{\pi(i)}|  d\pi \kl
2\pl \frac{1}{k} \summ_{j=1}^k x^*_j \pl .\]
\end{cor}

{\bf Proof:} Apply theorem \ref{KS} to $\al_{ij}\lel \cases{ x_i
&if $j\le \frac{n}{k}$ \cr 0& else}$. Let $m=[\frac{n}{k}]\ge 1$
be the smallest integer with $km\le n\le k(m+1)$, then the
rearrangement  $\al^*$ satisfies
\[ \frac{1}{n} \summ_{j=1}^n \al^*_j \lel \frac{1}{n} \kla [\frac{n}{k}](\summ_{i=1}^k x^*_i) + (n-k[\frac{n}{k}]) x^*_{k+1} \mer  \pl .\]
By the monotonicity of the $x^*_j$, we get \for \frac{1}{2}
\frac{1}{k} \summ_{i=1}^k x^*_i &\le& [\frac{n}{k}]\frac{1}{n}
\summ_{i=1}^k x^*_i \kl
\frac{1}{n} \kla [\frac{n}{k}](\summ_{i=1}^k x^*_i) + (n-k[\frac{n}{k}]) x^*_{k+1} \mer  \\
&\le& \frac{1}{k} \summ_{i=1}^k x^*_i  + \frac{1}{n}
(k([\frac{n}{k}]+1)-k[\frac{n}{k}]) \pl
\frac{1}{k} \summ_{i=1}^k x^*_i\\
&\le& (1+\frac{k}{n}) \pl \frac{1}{k} \summ_{i=1}^k x^*_i\pl. \\[-1.5cm]
\mel\qed

{\bf Proof of the lower estimate in Theorem \ref{main}:} As
usual, the trick for the lower estimate is an appropriate Abel
summation. Indeed, let $x_1\ge x_2\ge x_3\ge \cdots \ge x_n\ge 0$
be a non-increasing  sequence. It is well-known and easy to
check, that
\[ \noo x \rrm_X \lel \sup_{y\in C}  \summ_{i=1}^n x_iy_i \pl,\]
where $C$ consists of those vectors in the unit ball $B_{X^*}$
which are again  positive and non-increasing. In this case, we
have
\[ \summ_{i=1}^n x_iy_i \lel y_n \summ_{i=1}^n x_i + \summ_{k=1}^{n-1} (y_j-y_{j+1}) \summ_{i=1}^k x_i  \pl.\]
In other terms for an arbitrary vector $x$
 \begin{eqnarray}
 \noo x\rrm_X &=& \sup_{y\in C}
 y_n \summ_{i=1}^n x^*_i + \summ_{k=1}^{n-1} (y_k-y_{k+1}) \summ_{i=1}^k x^*_i  \pl.
 \end{eqnarray}
Now, we consider independent random variables $f_1,..,f_n$ and
the random vector
\[ x(\om)\lel (f_1(\om),\cdots,f_n(\om)) \pl.\]
Using $(1)$, we obtain \for
\lefteqn{\intt_\Om \noo \summ_{i=1}^n f_i(\om) e_i\rrm_X d\mu }\\
& & \lel  \intt_\Om \sup_{y\in C}
\kla y_n \summ_{i=1}^n x^*_i(\om) + \summ_{k=1}^{n-1} (y_k-y_{k+1}) \summ_{i=1}^k x^*_i(\om) \mer d\mu(\om) \\
 & & \gl  \sup_{y\in C}  y_n \intt_\Om (\summ_{i=1}^n x^*_i(\om)) d\mu(\om)  + \summ_{k=1}^{n-1} (y_k-y_{k+1}) \intt_\Om (\summ_{i=1}^k x^*_i(\om)) d\mu(\om) \\
 & & \lel  \sup_{y\in C} y_n \intt_\Om \noo \summ_{i=1}^n f_i e_i\rrm_1  d\mu  + \summ_{k=1}^{n-1} (y_k-y_{k+1}) \intt_\Om |\summ_{i=1}^k f_i e_i|_k  d\mu \pl .
\mel This reduces the problem to the investigation of norms $|\pl
|_k$. Let us define the increasing sequence
\[ s_j \lel (n\intt_{\frac{j-1}{n}}^{\frac{j}{n}} h^*(s) ds) \pl.\]
According to Corollary \ref{kk}, by independence and Fubini's
theorem, we deduce from the proof of Geiss's inequality
\cite[Theorem 3.4]{Gei} and monotonicity of $h^*$ \for
 \frac{2}{k} \pl \intt_\Om |\summ_{i=1}^k f_i e_i|_k  d\mu &\ge&
 \intt_\Om \intt_{\Pi_n} \sup_{i\le \frac{n}{k}} |f_{\pi(i)}(\om)| \pl d\pi d\mu
 \lel \intt_{\Pi_n} \intt_\Om \sup_{i\le \frac{n}{k}}
 |f_{\pi(i)}(\om)| \pl
  d\mu d\pi  \\
 &=& \intt_\Om \sup_{i\le \frac{n}{k}} |f_i(\om)|
 \pl
 d\mu  \gl
  \frac{1}{\sqrt{2}}  \pl [\frac{n}{k}] \intt_0^{\frac{1}{[\frac{n}{k}]} } h^*(s) ds \\
 &\ge& \pl \frac{1}{2\sqrt{2}} \pl \frac{n}{k} \intt_0^{\frac{k}{n}} h^*(s)
 ds\lel
 \pl \frac{1}{2\sqrt{2}} \pl \frac{1}{k}
 \summ_{j=1}^k
 (n\intt_{\frac{j-1}{n}}^{\frac{j}{n}} h^*(s) ds) \\
 &=& \pl \frac{1}{2\sqrt{2}} \pl \frac{1}{k}  \summ_{j=1}^k s_j \pl .
 \mel
Combining this with the previous estimate, we obtain again from
$(1)$ \for
 \intt_\Om \noo \summ_{i=1}^n f_i(\om) e_i\rrm_X d\mu &\ge& \frac{1}{4\sqrt{2}} \pl
 \sup_{y\in C} y_n\kla \summ_{j=1}^n s_j\mer +
 \summ_{k=1}^{n-1} (y_k-y_{k+1}) \summ_{i=1}^k s_i \\
 &=& \frac{1}{4\sqrt{2}} \pl \noo \summ_{i=1}^n s_i e_i\rrm_X \pl .
 \mel
Trivially, we have $\noo x\rrm_X\ge \noo x\rrm_\infty$ and
therefore Geiss' inequality \cite{Gei} concludes the proof
 \for
\kla \intt_\Om \noo \summ_i f_i e_i\rrm_X^p d\mu
\mer^{\frac{1}{p}} &\ge& \kla \intt_\Om \noo \summ_i f_i
e_i\rrm_\infty^p d\mu \mer^{\frac{1}{p}} \gl 2^{-\frac{1}{p}} \pl
(n\intt_0^{\frac{1}{n}} h^*(s)^p ds )^{\frac{1}{p}} \pl.\\[-1.5cm]
\mel\qed

\section{The upper bound}

In this section, we will prove the upper bound. Given  $n$
independent random variables $f_1,..,f_n$, we may use the fact
that for every monotone increasing function $g$
 \begin{eqnarray}
 \intt g(f) d\mu &=&  g(0) +\intt_0^\infty g'(s)
 \mu(|f|>s) ds \lel \intt_0^1 g(f^*(t)) dt \pl .
 \end{eqnarray}
Therefore, we can assume that $f_1,..,f_n$ are defined on
$[0,1]^n$ and non-increasing. First, we split all the $f_i$'s
into three parts. Let $b=h^*(\frac{1}{n})$ and
 \[ f_i^1\lel f_i1_{\{|f_i|>b\}} \quad ,\quad f_i^2 \lel
 (f_i-f_i^1)1_{[0,\frac1n)} \pl.  \]
We put $f_i^3=f_i-f_i^1-f_i^2$. The estimate for the first two
parts uses Rosenthal's \cite{Ro} inequality $(k\in \{1,2\})$
 \for
 \kla \intt \noo \summ_{i=1}^n f_i^k e_i\rrm_X^p d\mu \mer^\frac{1}{p}
 &\le&
 \kla \intt \kla \summ_{i=1}^n |f_i^k|\mer^p d\mu
 \mer^\frac{1}{p}\\
 &\le& c(p) \kla \intt \summ_{i=1}^n |f_i^k| d\mu +
 \kla \intt \sup_i |f_i^k|^p d\mu \mer^\frac{1}{p}
 \mer  \pl . \mel
We note that according to \cite{JSZ}, we have $c(p)\le
\frac{c_0p}{1+\ln p}$.  Let us observe that the two
 sets
 \for A^1 &=& \{(t,s_1,...,s_n) \p|\p f_i(s_i)>b \}
 \lel
 \{ (t,s_1,...,s_n) \p|\p |h(t,s_1,..,s_n)|>h^*(\frac1n)\} \pl , \\
 A^2 &=& \{(t,s_1,...,s_n) \p|\p \frac{i-1}{n} \le t <\frac in \Rightarrow
 s_i\le \frac{1}{n}  \}
 \mel
have measure less than $\frac1n$ and therefore
 \[  \frac{1}{n} \summ_{i=1}^n \noo f_i^k\rrm_1
 \lel  \intt_{A^k} |h|  \kl n^{\frac1p-1} \kla \intt_0^{\frac 1n}
 h^*(s)^p ds  \mer^{\frac1p}  \lel  \frac{1}{n} \kla
 n\intt_0^\frac{1}{n} h^*(s)^p   ds
 \mer^{\frac{1}{p}}   \pl.\]
Moreover, according to \cite[Theorem 3.4]{Gei}
 \for
 \kla \intt \sup_i |f_i^k|^p d\mu \mer^{\frac{1}{p}}
 &\le& \kla \intt \sup_i |f_i|^p d\mu
 \mer^\frac{1}{p}   \kl 2^{1-\frac1p} \pl \kla
 n\intt_0^\frac{1}{n} h^*(s)^p   ds
 \mer^{\frac{1}{p}}  \pl . \mel
Hence, by H\"olders inequality for $k\in \{1,2\}$
 \begin{eqnarray}
 \kla \intt \noo \summ_{i=1}^n f_i^ke_i\rrm_X^p d\mu \mer^\frac{1}{p}
 &\le&
 4c(p) \pl  \kla n\intt_0^\frac{1}{n} h^*(s)^p  ds
 \mer^{\frac{1}{p}} \pl .
 \end{eqnarray}
The estimate of the third  part uses the following proposition.

\begin{prop} \label{matri}
Let $X=(\rz^n,\noo \pl \rrm_X)$ be a symmetric sequence space,
$(\al_{ij})_{i,j=1,..,n}$ be a matrix and
$\al^*_1,...,\al^*_{n^2}$ be the non-increasing rearrangement of
the matrix, then
 \[ \kla n^{-n}
 \summ_{j_1,..,j_n=1}^n \noo \summ_{k=1}^n
 \al_{kj_k} e_k\rrm_X^p \mer^\frac{1}{p} \kl c_0 p
 \pl \noo \summ_{k=1}^n \al^*_{(k-1)n+1} e_k \rrm_X
 \pl .\]
\end{prop}

To conclude the proof of the upper estimate in Theorem
\ref{main},  we apply the Proposition to the matrix
\[ \al_{ij}\lel \sup_{\frac{j}{n}<s\le \frac{j+1}{n}} f^3_i(s) \pl.\]
$(\al_{in}=0$.) For fixed $\om=(s_1,..,s_n)$ with $\frac{j_i}{n}<
s_i\le \frac{j_i+1}{n}$, we get \for \noo \summ_{i=1}^n
f^3_i(s_i) e_i \rrm_X \kl \noo \summ_{i=1}^n \al_{ij_i} e_i
\rrm_X \pl . \mel The probability of the set where $s_i\in
(\frac{j_i}{n},\frac{j_i+1}{n}]$ is $\frac{1}{n}$ and for $s_i\le
\frac1n$ the variable $f_i^3(s_i)$ vanishes. Hence, we get
 \for
 \kla \intt \noo \summ_{i=1}^n f^3_i e_i\rrm_X^p \mer^{\frac{1}{p}}
 &\le&
  \kla n^{-n} \summ_{j_1,..,j_n=1}^n \noo \summ_{i=1}^n \al_{ij_i}
 e_i\rrm_X^p \mer^\frac{1}{p} \kl  \frac{c_0 p}{1+\ln p} \pl \noo \summ_{k=1}^n \al^*_{(k-1)n+1}
 e_k \rrm_X \pl . \mel
On the other hand, let us consider the new variables
 \[ \tilde{f}_i \lel \summ_{j=1}^{n}
 \al_{ij} 1_{[\frac{j-1}{n},\frac{j}{n})} \pl.\] By
definition $\tilde{f}_i\le f_i$ and therefore
 \[ \tilde{h}(t,s_1,..,s_n) \pl:=\pl  \summ_{i=1}^n
 1_{[\frac{i-1}{n},\frac{i}{n}]} \tilde{f}_i  \kl h
 \pl .\]
We observe that the non-increasing rearrangement of $\tilde{h}$
is the same as the non-increasing rearrangement of the matrix
$\al$. For simplicity, let us assume that the values $\al_{ij}$'s
are all different from each other. These values appear in
$\tilde{h}$ on the disjoint sets
\[ A_{ij} \lel \left\{(t,s_1,..,s_n) \mitt
\frac{i-1}{n}\le t < \frac{i}{n} \pl ,\pl \frac{j-1}{n}\le s_i <
\frac{j}{n} \right\} \] of measure $\frac{1}{n^2}$. This implies
for all $j=1,..,n$
\[ \al^*_{jn+1} \lel \tilde{h}^*(\frac{j}{n})  \pl .\]
Together with
 \[ \tilde{h}^*(\frac{j}{n}) \kl n\intt_{\frac{j-1}{n}}^{\frac{j}{n}}
 \tilde{h}^*(s)ds \kl n\intt_{\frac{j-1}{n}}^{\frac{j}{n}}
 h^*(s)ds  \pl ,\]
we deduce  from $\al^*_1\le b=h^*(\frac{1}{n})$ that

\begin{samepage} \for \noo \summ_{k=1}^n \al^*_{(k-1)n+1} e_k\rrm_X &\le&
 \al^*_1 +
 \noo \summ_{k=2}^n h^*(\frac{k-1}{n}) e_k\rrm_X \\
 &\le& h^*(\frac{1}{n}) + \noo \summ_{k=2}^n
 (n\intt_{\frac{k-2}{n}}^{\frac{k-1}{n}} h^*(s)ds) e_k\rrm_X \\
 &\le& n\intt_0^{\frac{1}{n}} h^*(s) ds   + \noo \summ_{k=1}^n
 (n\intt_{\frac{k-1}{n}}^{\frac{k}{n}} h^*(s)ds) e_k\rrm_X \\
 &\le& (n\intt_0^{\frac{1}{n}} h^*(s)^p ds)^{\frac{1}{p}}    + \noo
 \summ_{k=1}^n (n\intt_{\frac{k-1}{n}}^{\frac{k}{n}} h^*(s)ds)
 e_k\rrm_X \pl.
 \mel \end{samepage}

This concludes the proof of (Proposition 2.1 $\Rightarrow$ upper
estimate in Theorem \ref{main}). The proof of Proposition
\ref{matri} relies on the combinatorial estimate \ref{comb} and
the following elementary observation.

\begin{lemma} \label{elel}
Let $x=(x_1,..,x_n)$ be a positive non increasing sequence $y\in
\rz^n$ built by repetitions of the coordinates  in $x$, i.e. for
$y_j\in\{x_1,..,x_n\}$ for all $1\le j\le n$. If
\[ \beta_i \pl:=\pl {\rm card}\{ j \p|\p  y_j\lel x_i \},\]
then
\[ \noo y\rrm \kl 2\pl \max\{1,\sup_r \frac{1}{r} \summ_{i=1}^r \beta_i
\}\pl \noo x\rrm \pl .\]
\end{lemma}

{\bf Proof:} We assume
\[  \summ_{i=1}^{r} \beta_i \le t r\]
for all $r$ and $t\gl 1$. Let $1\le k\le n$. In order to
calculate $y^*_k$, we choose $r$ such that
\[ \summ_{i=1}^{r-1} \beta_i \pl
<\pl k \kl \summ_{i=1}^r \beta_i
  \pl .\]
By the assumption on $y$, this means $y^*_k=x_r$. By assumption
on $t$, we have
\[ k\le \summ_{i=0}^r \beta_i \kl tr \pl .\]
In other terms $y^*_k\kl x_r \kl x_{[\frac{k}{t}]}$. Here we use
the convention $[r]=1$ for $r<1$. In order to estimate the norm
of $T_t(x)=x_{\frac{k}{t}}$, we fix a nonincreasing sequence $z$
in the unit ball of $X^*$. We use Abel summation for
$\tilde{x}=T_t(x)$ \for
 \langle \tilde{x},z\rangle &=& z_n(\summ_{i=1}^{n} \tilde{x}_i) +
 \summ_{k=1}^{n-1} (z_k-z_{k+1}) \summ_{j=1}^k \tilde{x}_j \pl.
 \mel
We fix $1\le k\le n$, then
 \for \summ_{j=1}^k
  x_{[\frac{j}{t}]}  &\le&  tx_1 + \summ_{1\le i\le
 \frac{k}{t}} tx_i \kl 2t \summ_{i\le k} x_i \pl.
 \mel
Hence, we deduce from $z_k-z_{k+1}\ge 0$ that
 \for
 \langle \tilde{x},z\rangle &\le& 2t \pl \kla
 z_n(\summ_{i=1}^{n} x_i) + \summ_{k=1}^{n-1}
 (z_k-z_{k+1}) \summ_{j=1}^k x_j \mer \kl  2t \langle
 x,z\rangle \kl 2t \noo z\rrm_{X^*}\pl \noo x\rrm_X
 \pl .\mel
This shows
 \for \noo y\rrm_X \kl \noo T_t(x)\rrm
 \lel \sup_z \langle T_t(x),z\rangle
 \kl 2t \pl \noo x\rrm \pl .\\[-1.5cm]
 \mel\qed

{\bf Proof of Proposition \ref{matri}:} The combinatorial
estimate will be shown in the next paragraph. We will show how it
can be used to prove the Proposition. Using $(1)$, we can again
assume $\al_{i1}\gl \al_{i2} \gl \cdots \gl \al_{in}$ for all
$i=1,..,n$. Moreover, by perturbation, we can assume that all the
$\al_{ij}$'s are different. For fixed $i,k$ we denote by $I_{ik}$
the interval in $[1,n]$ satisfying
 \[ j\in I_{ik} \quad \Leftrightarrow \quad \al^*_{kn+1} \pl < \pl \al_{ij} \kl
 \al^*_{(k-1)n+1} \pl .\]
We define
 \[ \mu_{ik} \pl:=\pl \frac{1}{n} {\rm card}(I_{ik}) \pl.\]
Using the convention $\al+{n^2+1}=0$, we have
 \[ \summ_{k=1}^n \mu_{ik} \lel 1 \pl .\]
On the other hand, there are $n$ different values between
$\al^*_{kn+1}$ and $\al^*_{(k-1)n+1}$, hence
 \[ \summ_{i=1}^n \mu_{ik} \lel 1 \pl .\]
Thus $\mu$ is a doubly stochastic matrix. For $j\in I_{ik}$, we
can replace $\al_{ij}$ by the bigger value $\al^*_{(k-1)n+1}$.
Let us denote this modified matrix by $\tilde{\al}$. Let us
observe that the random variable
 \[ \tilde{f}_i \lel \summ_{f=1}^n \tilde{\al}_{ij} 1_{\{j\}} \]
on the probability space $\{1,..,n\}$ with
$\la(\{j\})=\frac{1}{n}$ satisfies
 \[ \la\kla \tilde{f}_i\lel \al^*_{(k-1)n+1)} \mer \lel \mu_{ik} \pl.\]
Hence $f_i$ has the same distribution (and hence the same
non-increasing rearrangement) as
\[ g_i \lel \summ_{k=1}^n \al^*_{(k-1)n+1} 1_{\{k\}} \quad \mbox{
with respect to the measure } \mu_i(\{k\})\lel \mu_{ik} \pl .\]
Therefore, the assertion follows from an estimate of
\[ \kla \intt_{\{1,..,n\}^n} \noo \summ_{i=1}^n g_i(j_i) e_i \rrm_X^p
d\mu_1(j_1)\cdots d\mu_n(j_n) \mer^{\frac{1}{p}} \pl .\] For an
individual element $\om=(j_1,..,j_n)$, we observe that
$(g_1(j_1),\cdots,g_n(j_n))$ only takes the values in the set
$\al^*_{(k-1)n+1}$ for $k=1,..,n$. Moreover, let  $\beta_k$ be
the cardinality  of this occurrence, then
\[ \beta_k \lel \summ_{i=1}^n \delta_{kj_i}  \pl .\]
Hence, we get
 \[ \frac{1}{r} \summ_{k=1}^r \beta_k \lel \frac{1}{r} \summ_{i=1}^n
 1_{[1,r]}(j_i) \pl .\]
According to Lemma \ref{elel}, we have
 \[ \noo \summ_{i=1}^n g_i(j_i) e_i\rrm_X \kl 2 \max\{1,\sup_r
 \frac{1}{r} \summ_{i=1}^n 1_{[1,r]}(j_i)\} \pl \noo \summ_{k=1}^n
 \al^*_{(k-1)n+1} e_k\rrm_X  \pl . \]
Theorem \ref{comb} yields the assertion.\qed

\begin{rem} {\rm The typical example for the theorem is the (quasi-) norm
\[ \noo x\rrm_{1,\infty} \lel \sup_k k \pl x^*_k \pl .\]
If we consider $\al_{ij}=\frac{1}{j}$ and thus
$\al^*_{(k-1)n+1}=\frac{1}{k}$, then the norm of a repetition
$y_k(j_1,...,j_n)=j_k^{-1}$ is exactly
\[ \noo \summ_{k=1}^n y_ke_k \rrm_{1,\infty} \lel \sup_r \frac{1}{r} \summ_{i=1}^r \beta_i \pl
.\] Moreover, the estimate
\[ \noo T_t:\ell_{1,\infty}\to \ell_{1,\infty}\rrm \kl 2\max\{1,t\}\]
is still valid. In this case, we see that Proposition \ref{matri}
for $X=\ell_{1,\infty}^n$ is equivalent to the combinatorial
estimate in Theorem \ref{comb}.}
\end{rem}

\begin{rem}{\rm  Let $f$ be an increasing function and consider
the Lorentz space
 \[ \noo x\rrm_{f,w} \lel \kla \summ_{k=1}^n \kla
 \frac{f(k)}{k^{\frac 1w}} x^*_k\mer^w \mer^{\frac1w} \quad
 \mbox{and} \quad \noo x\rrm_{f,\infty} \lel \sup_k f(k)x_k^* \pl .
 \]
In this context Hardy's inequality reads as follows.  If
 \for
 \frac{f(n)}{f(k)} &\le&  c_s \kla \frac{n}{k} \mer^{\frac1s}
 \mel

and $q<s$, we may find $\gamma>0$ such that $p(\frac1s+\gamma)<1$.
Then, for some constant $c(\gamma)$ we have
 \[ \kla \summ_{k=1}^n [\frac{f(n)}{f(k)} k^{-\gamma}]^q
 \mer^{\frac1q} \kl c_s c(\gamma)  n^{\frac1q-\gamma} \pl .\]
Following \cite[2.1.7, p=75]{Pie}, we get
 \[ \kla \summ_{k=1}^n {x_k^*}^q \mer^{\frac1q} \kl
 c_s c(\gamma) \frac{n^{\frac1q}}{f(n)n^\gamma} \kla \summ_{k=1}^n
 [f(k)k^{\gamma-\frac1w} x_k^*]^w \mer^{\frac1w} \pl .\]
The same calculation as in \cite[2.1.7]{Pie} then yields the
Hardy inequality
 \[ \noo x\rrm_{l_{f,w}} \kl
 \noo \summ_{k} \kla \frac{1}{k} \summ_{j=1}^k {x_j^*}^q \mer^{\frac1q} e_k \rrm_{l_{f,w}}
 \kl c'(\gamma)c_s
 \noo x\rrm_{l_{f,w}} \pl .\]
Moreover, of we assume in addition $q\le w$, we may combine this
argument with  the triangle inequality in $\ell_{\frac{w}{q}}$
and deduce the $q$-convexity of $l_{f,w}$, i.e.
 \begin{eqnarray}
 \noo \kla \summ_{j} |x_j|^q \mer^{\frac1q} \rrm_{l_{f,w}} &\le&
 c_s c'(\gamma) \kla \summ_j \noo x_j\rrm_{l_{f,w}}^q \mer^{\frac1q} \pl .
 \end{eqnarray}
Therefore
 \[ \noo| x|\rrm \lel \noo |x|^{\frac1q} \rrm_{l_{f,w}}^{q} \]
is equivalent to a norm.  Given $p\ge q$, we may then apply
Theorem \ref{main} to $(\rz^n, \noo | \pl |\rrm)$ and $|f_i|^q$
and obtain the lower estimate (with $(c_s
c'(\gamma)1(1+\sqrt{2}))^{\frac1q}$) and the  upper estimate with
$(c_s c'(\gamma) c(\frac{p}{q}))^{\frac1q}$.

Now, let us consider more generally a symmetric quasi-norm
satisfying the Hardy inequality
 \[ \noo x\rrm_X \kl \noo \summ_{k=1}^n (\frac{1}{k} \summ_{j=1}^k
 {x_j^*}^q)^{\frac1q} e_k \rrm_X \kl c(q) \noo x\rrm_X \pl .\]
Then, we can easily modify Lemma \ref{elel} and deduce that for
non-increasing $x$ and $y$ with
 \[  \beta_i \pl:=\pl {\rm card}\{ j \p|\p  y_j\lel x_i \} \]
we have
 \for
 \noo y\rrm_X \kl c(q)^{\frac1q} \pl 2^{\frac1q} \pl
  \max\{1,\sup_r \frac{1}{r} \summ_{i=1}^r \beta_i \}^{\frac1q} \pl \noo x\rrm_X \pl .
   \mel
Therefore Proposition 2.1. is still valid in this setting (using
the combinatorial estimate for $\frac{p}{q}$ which is big when
$q$ is small.) If $X$ is quasi-normed, there is an equivalent norm
$\noo | \pll | \rrm$ satisfying $\noo| x+y|\rrm^r\kl \noo |
x|\rrm^r+ \noo |y|\rrm^r$. Then the argument using Rosenthal's
inequality easily works for $p\ge r$ and the constant
$c(\frac{p}{r})^{\frac1r}$. Therefore the upper estimate holds
with $c_0 \max\{(\frac{p}{r(1+\ln \frac{p}{r})})^{\frac1r},
(\frac{pc_s}{q(1+\ln \frac{p}{q})})^{\frac1q}\}$ provided Hardy's
inequality is available. At the time of this writing it is not
clear whether the lower estimate still holds under these
assumptions.}

\end{rem}

\section{The combinatorial estimate}
The combinatorial estimate is based on a tail estimate for the
variables
\[ h_r(j_1,..,j_n) \lel \summ_{i=1}^n 1_{[1,r]}(j_i) \pl.\]
with respect to the product probability measure on $\{1,..,n\}^n$
defined by
\[ P_\mu\{(j_1,..,j_n)\} \lel  \prod_{k=1}^n \mu_{kj_k}\pl .\]
Note again, that
 \for
  \int_{\{1,...,n\}^n} \noo \summ_{l=1}^n
 j_l^{-1}e_l\rrm_{1,\infty}^p d P_{\mu}(j_1,...,j_n)
  \lel \int_{\{1,...,n\}^n} \sup_{r}
   (\frac{h_r}{r})^p \pl dP_{\mu} \pl .
  \mel
Therefore, the combinatorial estimate is a special case of our
main results for the weak-$\ell_1$ `norm'. Here
$\mu=(\mu_{ij})_{ij}$ is assumed to be a doubly stochastic matrix
and $P_\mu$ denotes the product probability measure on
$\{1,..,n\}^n$ defined by
\[ P_\mu\{(j_1,..,j_n)\} \lel  \prod_{k=1}^n \mu_{kj_k}\pl .\]
Let us denote by $\C\subset \rz^n$ the set of all doubly
stochastic matrices. The following lemma provides the key
estimate using Birkhoff's theorem on doubly stochastic matrices.
(Although the intuition for this estimate comes from the
non-extremal matrix $\mu_{ij}=\frac1n$.)

\begin{lemma} \label{Herz}
Let $\mu\in C$. Then
 \[ P_\mu(h_r=j) \kl \summ_{\stackrel {k=0}{
 (n-j)\ge (j-k), k\le j} }^{r} \kla \frac kj \mer^j
 \pl \kla \frac{j-k}{n-j}\mer^{n-j} \kla {r \atop k}
 \mer \pl \kla {n-j \atop j-k} \mer \pl \]
holds for all $1\le j\le r \le  n$ and $j<n$. For $j=r=n$,
$P_\mu(h_n=n)=1$.
\end{lemma}

{\bf Proof:} The equality $h_r=j$ holds if and only if there is a
set $B\subset \{1,..,n\}$ of cardinality $r$ such that $j_i\le r$
for $i\in B$  and $j_i>r$ for $i\notin B$. Using the
geometric/arithmetic mean inequality, we deduce
 \for P_\mu(h_r=j) &=& \summ_{|B|=j} \prod_{i\in B}
 \kla \summ_{s\le r} \mu_{is}
 \mer \pl \prod_{i\notin B} \kla \summ_{s>r} \mu_{is} \mer  \\
 &\le& \summ_{|B|=j} \kla \frac1j \summ_{i\in B,
 s\le r} \mu_{is} \mer^j  \pl \kla \frac{1}{n-j}
 \summ_{i\notin B, s>r} \mu_{is} \mer^{n-j} \pl .
 \mel
Since $\mu$ is doubly stochastic, we can simplify the second term
 \for \summ_{i\notin B, s>r}
 \mu_{is} &=& \summ_{i\notin B}
 \kla 1- \summ_{s\le r} \mu_{is} \mer \lel
  \summ_{i=1}^n \kla 1- \summ_{s\le r} \mu_{is} \mer -\summ_{i\in B}
 \kla 1- \summ_{s\le r} \mu_{is} \mer \\
 &=& n - \summ_{s\le r} \summ_{i=1}^n \mu_{is} - j+\summ_{i\in B,s\le r}
 \mu_{is} \\
 &=& n-r-j+ \summ_{i\in B,s\le r} \mu_{is}  \pl .
 \mel
We are lead to  the function
 \[ f(x)\lel
 \kla \frac xj\mer^j \pl \kla 1-\frac{r-x}{n-j} \mer^{n-j}
 \pl .\]
It is elementary to check that
 \for f''(x) &=&
 \frac{j-1}{j}\kla \frac xj\mer^{j-2} \kla
 1-\frac{r-x}{n-j}\mer^{n-j} + 2 \kla \frac
 xj\mer^{j-1}
 \kla 1-\frac{r-x}{n-j}\mer^{n-j-1} \\
 & & \pll + \frac{n-j-1}{n-j} \kla
 1-\frac{r-x}{n-j}\mer^{n-j-2}\kla \frac xj\mer^j
 \mel
which is positive on the interval $I=[\max\{0,(r+j-n)\},\infty)$.
Indeed, we note $x\ge (r+j-n)$ iff $1-\frac{r-x}{n-j}\ge 0$. On
the convex set  $\C$, we consider the linear functional
\[ L(\mu) \lel \summ_{i\in B} \summ_{s\le r} \mu_{is} \gl 0 \pl .\]
Using again
 \[ 1 -\frac{r-L(\mu)}{n-j} \lel \frac{1}{n-j} \kla n-r-j +\summ_{i\in B,s\le r} \mu_{is} \mer
 \lel  \frac{1}{n-j} \summ_{i\notin B} \summ_{s>r}
 \mu_{is} \gl 0 ,\]
we deduce that $L(\mu)\in
 [\max\{0,(r+j-n)\},\infty)$. Fixing a subset $B$, we deduce that
 \[ h_B(\mu) \lel f(L(\mu)) \lel f(\summ_{i\in B,s\le r} \mu_{is}) \]
is a convex function on $\C$. The case $j=n$ is excluded by
assumption.  In particular,
 \[ h(\mu) \lel \summ_{|B|=j} f(\summ_{i\in B, s\le r} \mu_{is}) \]
is convex and attains its maximum on an extreme point. According
to \cite{Bi}, the extreme points in $\C$ are the permutation
matrices. For any permutation $\pi:\{1,..,n\}\to \{1,..,n\}$ we
have
 \[ \summ_{i\in B} \summ_{s\le r} \delta_{\pi(i)s}
 \lel \summ_{i \in \pi^{-1}(B)} \summ_{s\le r}
 \delta_{is} \lel |\pi^{-1}(B)\cap \{1,..,r\}| \pl .
 \]
Since the map $B\mapsto \pi^{-1}(B)$ yields a bijection on the
subsets of $\{1,..,n\}$  of cardinality $j$, it is sufficient to
consider the trivial permutation $\pi(i)=i$ and thus the identity
matrix $\delta=(\delta_{ij})$. We define $A_r=\{1,..,r\}$ and get
  \for h(\delta) &=&
 \summ_{|B|=j} \kla \frac{1}{j} \summ_{i\in B}
 \summ_{s\le r} \delta_{is}\mer^j \kla \frac{1}{n-j}
 \summ_{i \notin B} \summ_{s>r}
 \delta_{is} \mer^{n-j} \\
 &=&
  \summ_{|B|=j} \kla \frac{|A_r\cap B|}{j}\mer^{j} \pl
 \kla \frac{|B^c\cap A_r^c|}{n-j}\mer^{n-j} \\
 &=& \summ_{k=0}^r \summ_{\stackrel{C\subset A_r, |C|=k,}{D\subset A_r^c,
 |D|=j-k}}\kla \frac kj\mer^j \kla \frac{j-k}{n-j} \mer^{n-j} \\
 &=&
 \summ_{k=0, (n-j)\ge (j-k),k\le j }^{r}
 \kla \frac kj \mer^j \pl \kla
 \frac{j-k}{n-j}\mer^{n-j} \kla { r \atop k} \mer \pl \kla {n-j \atop j-k}
 \mer \pl .
 \mel
Since $P_\mu(h_r=j)$ is  majorized by $h(\delta)$, this concludes
the proof.\qed

\begin{lemma} Let $t\gl e^2$ and $r<n$, then
\[ P_\mu(h_r\gl tr) \kl 2  \pl \kla \frac{e^3}{t}\mer^{tr} \pl.\]
\end{lemma}

{\bf Proof:} According to Lemma \ref{Herz}, we have \for
 P_\mu(h_r\ge tr) &=& \summ_{n>j\gl tr} P_\mu(h_r=j) \\
 &\le& \summ_{n>j\ge tr} \summ_{k=0, (n-j)\ge (r-k) }^{\min(r,j)}
 \kla \frac kj \mer^j \pl \kla
 \frac{j-k}{n-j}\mer^{n-j} \kla { r \atop k} \mer \pl \kla {n-j \atop j-k}
 \mer \\
 &=&
 \summ_{k=0}^r \kla { r \atop k} \mer \summ_{n>j\ge tr,k\le j, (n-j)\ge
 (r-k)} \kla \frac kj \mer^j \pl \kla \frac{j-k}{n-j}\mer^{n-j} \pl  \kla
 {n-j \atop j-k} \mer  \pl.\mel
Using Stirling's formula ($m!=m^me^{-m}\sqrt{2\pi
m}e^{\frac{\theta_m}{12}}$, $0\le \theta_m\le 1$) we deduce (with
$e^{\frac{1}{12}}\le \sqrt{2\pi}$) that
 \[ \kla { n \atop k} \mer \le  \frac{n^n}{(n-k)^{n-k}
 k^k} \pl ,\]
we get for $n-j\ge j-k$ \for
 \lefteqn{\kla
 \frac{j-k}{n-j}\mer^{n-j} \pl  \kla {n-j \atop j-k} \mer }\\
 & & \kl  \kla \frac{j-k}{n-j}\mer^{n-j}
 \frac{(n-j)^{n-j}}{(n-j-(j-k))^{(n-j-(j-k))} (j-k)^{(j-k)} } \\
 & & \lel   \kla
 \frac{j-k}{n-j-(j-k)} \mer^{n-j-(j-k)} \\
 & & \lel \kla 1+
 \frac{(j-k)-(n-j)+(j-k)}{n-j-(j-k)}\mer^{n-j-(j-k)} \\
 & & \kl \max\{1,\exp(2(j-k)-n+j)\} \kl \exp(2(j-k)) \pl . \mel
Therefore, we have
 \[   P_\mu(h_r\ge tr) \kl \summ_{k=0}^r \kla { r \atop k} \mer \exp(-2k)
 \summ_{j\ge tr} \kla \frac{ke^2}{j} \mer^j \pl .\]
For fixed $k$, we consider $f(x)=\exp(x[2+\ln k-\ln x])$ which
satisfies $f'(x)=f(x)[1+\ln k-\ln x]$. On the interval
$[ke,\infty)$ the function $f$ is decreasing. For $x\ge ke$, we
consider
\[ g(x) \lel  -\frac{f(x)}{\ln x-1-\ln k} \]
and observe
 \[ g'(x)\lel f(x)\kla 1
 +\frac{1}{x[\ln x-1-\ln k]^2}\mer \gl f(x) \pl .\]
Let $j_0$ be  such that $tr\le j_0 < tr+1$, then we deduce from
the monotonicity of $f$ and with $j_0\ge tr \ge e^2r\ge e^2k$ that
 \for \summ_{j\ge j_0} \kla
 \frac{ke^2}{j} \mer^j &\le& \kla \frac{ke^2}{j_0}
 \mer^{j_0}+ \intt_{j_0}^\infty
 f(x) \pl dx \kl
  \kla \frac{ke^2}{j_0} \mer^{j_0} +
 \intt_{j_0}^\infty g'(x) \pl dx \\
 &=& \kla \frac{ke^2}{j_0} \mer^{j_0} + \frac{1}{\ln j_0-\ln ke} \pl
 \kla \frac{ke^2}{j_0} \mer^{j_0} \\
 &\le & \kla \frac{ke^2}{j_0} \mer^{j_0} + \kla \frac{ke^2}{j_0}
 \mer^{j_0}  \lel 2 \kla \frac{ke^2}{j_0} \mer^{j_0}
 \kl
   2
 \kla \frac{ke^2}{tr} \mer^{tr} \pl .\mel
Since  $t\ge 2$, we deduce
 \for
 P_\mu(h_r\ge tr) &\le& \summ_{k=0}^r \kla {k
 \atop r}\mer \pl
 \exp(-2k) \summ_{j\ge tr} \kla \frac{ke^2}{j} \mer^j
 \kl
    \summ_{k=0}^r \kla {k \atop r}\mer \pl \exp(-2k) \pl
 2\kla \frac{ke^2}{tr} \mer^{tr}\\
 &\le& 2^{2r} \kla \frac{e^2}{t} \mer^{tr}  \kl
  \kla \frac{e^3}{t} \mer^{tr} \pl .\\[-1.5cm]
 \mel\qed

The next calculation provides the $\frac{p}{\ln p}$ term.

\begin{lemma}\label{cc1} Let $b\ge 1$ and $a\ge
\max\{\frac{e^{e-1}b}{2},4b^2\}$. Then
 \[ e^{-2a-1} \kla \frac{a}{1+\ln a}
    \mer^{a+1} \kl \int_b^\infty
    t^a \kla \frac{b}{t}\mer^{t} dt \kl (a+1)
    \kla \frac{2a}{1+\ln(a)}\mer^a  \pl .\]
\end{lemma}

{\bf Proof:} The derivative of the function $f(t)=t^a\kla
\frac{b}{t}\mer^t$ is given by $f'(t)=f(t)(\frac{a}{t}+\ln b-\ln
te)$ and thus $f$ has a  unique maximum for $t_0$ satisfying
 \[ t_0 (\ln t_0e-\ln b) \lel a \pl .\]
Let us denote by $M=\sup f(t)$.  Then we have
 \[ M \lel f(t_0) \lel t_0^ae^{-t_0(\ln t_0-\ln b)}
 \lel e^{t_0}t_0^ae^{-t_0(\ln t_0e-\ln b)} \lel
 e^{t_0-a}t_0^a \kl t_0^a \pl .\]
Note that $g(t)=\kla\frac{b}{t}\mer^t$ is decreasing on
$(b,\infty)$ and thus
 \for
  \intt_b^\infty t^a \kla \frac{b}{t}\mer^{t} dt
 &\gl&  \intt_{\frac{t_0}{2}}^{t_0} t^a \kla \frac{b}{t} \mer^t
 dt \gl
    \kla \frac{t_0}{2}\mer^{a+1} \kla \frac{b}{t_0}
 \mer^{t_0} \lel 2^{-(a+1)} t_0 M \gl 2^{-(a+1)} e^{-a} t_0^{a+1}    \pl.
 \mel
Hence the lower estimate follows from a lower estimate of $t_0$.
Indeed, the function $l(t)=t(\ln te-\ln b)$ is increasing on
$(b,\infty)$ and it is easily checked that $b,a\ge 1$ implies
$l(\frac{a}{\ln ea})\le a$. Hence, $t_0\ge \frac{a}{\ln(ea)}$ and
also $a\ge 4b^2$ implies $\frac{a}{\ln ea}\ge 2b$. We deduce
 \[   \intt_b^\infty t^a \kla \frac{b}{t}\mer^{t} dt
    \gl e^{-a}2^{-(a+1)}  \kla \frac{a}{1+\ln a}
    \mer^{a+1} \pl .\]
For the upper estimate of  $t_0$ we note that $s\ge e$ implies
$\frac{\ln s}{s}\le e^{-1}\le \frac12$. Therefore, we deduce for
$\gamma=\frac{2e}{b}$  that
 \for
  l(\frac{2a}{\ln(\gamma a)})
  &=& \frac{2a}{\ln (\gamma a)}
  [\ln(\gamma a)-\ln
  \ln(\gamma a)] \lel  2a[1-\frac{\ln \ln (\gamma
  a)}{\ln(\gamma a)}]\gl a \pl .
  \mel
By our assumptions $\frac{a}{\ln(\gamma a)}\kl \frac{2a}{1+\ln
a}$. This yields $M\le t_0^a\kl \kla \frac{2a}{1+\ln(a)}\mer^a$.
For the upper estimate of the integral, we consider $b(t)=
\frac{a}{t}+\ln b-\ln(te)$ and $h(t)\lel \frac{f(t)}{b(t)}$. Then
the derivative of $h$ satisfies
 \[ h'(t)\lel f(t)[1-\frac{b'(t)}{b(t)^2}]
 \lel f(t)[1+\frac{
 \frac{a}{t^2}+\frac{1}{t}}{b(t)^2}] \gl f(t) \pl
 .\]
Since $a(t)$ is negative for $t\ge a$, we deduce from
$\frac{a}{b}\ge e$ that
 \for
 \int_b^\infty
    t^a \kla \frac{b}{t}\mer^{t} dt
    &=& \intt_b^a t^a \kla \frac{b}{t}\mer^{t} dt
    + \intt_a^\infty t^a \kla \frac{b}{t}\mer^{t} dt
    \\
   &\le& aM + \intt_a^\infty h'(t)  dt \\
   &\le& a M + \frac{f(a)}{\ln(ae)-\ln b-1} \kl
    (a+1)M \pl .
 \mel
The assertion is proved. \qed

 The next lemma is elementary. The proof
 uses $\frac{b}{t}\le e^{-1}$ for $t\ge eb$ and is
easier than the proof of Lemma \ref{cc1}, we omit the details.

\begin{lemma} \label{calc} Let $b\ge 1$ and  $a>0$, $d\ge be$, then
\[ \intt_d^\infty \kla \frac bt\mer^{tr} t^a \pl dt \kl \frac{2}{r} \pl
\exp(-dr) \pl \cases{d^a & if $2a\le rd$ \cr 2\kla \frac{2a}{r}
\mer^a &if $rd\le 2a$ } \pl . \]
\end{lemma}

The proof of the combinatorial result  is now rather a matter of
calculation. \hz

{\bf Proof of Theorem \ref{comb}:} Let $d=e^4\gl 2$, and choose
$r_0$ such that $r_0d\le 2(p-1)\le (r_0+1)d$. We also use $b=e^3$
and assume $p\ge 2 \max\{e^{e+2},4e^6\}=8e^6$ (which implies
$r_0+1\ge 8e^2$). Then, we deduce from Lemma \ref{cc1} and Lemma
\ref{calc}
 \for
 \lefteqn{\intt \kla \sup_{r<n} \frac{h_r}{r} \mer^p dP_\mu }\\
 & &\lel p \intt_0^\infty t^{p-1} P_\mu (\bigcup_r \{h_r> rt\}) \pl dt \\
 & &\kl p\intt_0^d t^{p-1} dt+\intt_d^\infty  t^{p-1} \summ_{r=1}^\infty
 P_\mu(h_r\gl rt)\pl  dt \\
 & & \kl d^p + \summ_{r=1}^\infty  p \intt_d^\infty \kla \frac{e^3}{t}
 \mer^{tr} \pl t^{p-1} \pl dt \\
 & & \kl d^p +  p \kla \summ_{r=1}^{r_0}   \intt_d^\infty \kla
 \frac{e^3}{t} \mer^{tr} \pl t^{p-1} \pl dt +
 \summ_{r>r_0}^{\infty}  \intt_d^\infty \kla
 \frac{e^3}{t} \mer^{tr} \pl t^{p-1} \pl dt \mer  \\
 & & \kl d^p + p \kla \summ_{r=1}^{r_0}  p \kla \frac{2(p-1)}{1+\ln(p-1)} \mer^{p-1}  +
 \summ_{r>r_0} d^{p-1} \frac{2}{r} \exp(-dr) \mer\\
 & & \kl d^p + r_0 p^2  \kla \frac{2(p-1)}{1+\ln(p-1)} \mer^{p-1}
 + pd^{p-1} (\summ_{r\ge 1} \exp(-dr))\\
 & & \kl d^p +  p^3  \kla \frac{4p}{1+\ln p} \mer^{p-1}
 + 2pd^{p-1} \pl .
 \mel
Using the triangle inequality in $L_p$ and $\ell_p^3$ and
$x^{\frac1x}\le e^{\frac1e}$, we conclude the proof
 \for
 \kla \intt \kla \sup_{r\le n} \frac{h_r}{r} \mer^p
 dP_\mu\mer^{\frac1p} &\le& \kla \intt \frac{h_n}{n}^p
 dP_{\mu} \mer^{\frac 1p} +
  \kla \intt \kla \sup_{r<n} \frac{h_r}{r} \mer^p
 dP_\mu\mer^{\frac1p} \\
 &\le& 1 + d + p^{\frac{3}{p}}
 \kla \frac{4p}{1+\ln p} \mer^{\frac{p-1}{p}} +
 (2p)^{\frac1p} d^{\frac{p-1}{p}} \\
 &\le& 1+ d+ e\frac{4p}{1+\ln p} + ed  \kl 2e^5+
 e\frac{4p}{1+\ln p} \pl . \\[-1.5cm]
 \mel\qed

\begin{rem} For $p=1$, we can use the first part of Lemma \ref{calc} and get the `concrete
estimate'
\[ \intt \sup_r \frac{h_r}{r} dP_\mu \kl 2+e^4 \pl .\]
\end{rem}

\begin{rem} {\rm As a standard  application, we obtain a fairly good tail
estimate. Assuming $\noo \summ_{k=1}^n \al^*_{(k-1)n+1}e_k
\rrm_X=1$, we have
 \[ Prob \kla  \noo \summ_{k=1}^n \al_{kj_k}e_k \rrm_X>
 t\mer  \kl C \exp(-\frac{t\ln t}{C}) \]
for some universal constant $C$. As usual this is obtained from
Chebychev's inequality and choosing $p$  optimal. We obtain a
similar behaviour for general independent  functions bounded by
$1$ and such that $\noo \summ_{i=1}^n
(n\intt_{\frac{i-1}{n}}^{\frac
 in} h^*(s) ds ) \pl e_i \rrm_X \le 1$. Of course, this behaviour
 is a reformulation of our main result.}
\end{rem}

\begin{exam} Let $\mu$ be the standard matrix $\mu_{ij}=\frac1n$ and $x=e_1$ the first unit vector and $\al_{ij}=
\cases{1 & if $j=1$ \cr 0 & else}$. Then
 \[
 \kla n^n \summ_{j_1,...,j_n=1}^n \noo \summ_{k=1}^n
 x(j_k) e_k\rrm_{\ell_1^n}^p \mer^{\frac1p} \lel
 \kla n^n \summ_{j_1,...,j_n=1}^n \noo \summ_{k=1}^n
 \al_{k,j_k} e_k\rrm_{\ell_1^n}^p \mer^{\frac1p}
 \lel \kla \intt h_1^p dP_{\mu} \mer^{\frac1p}  \]
Moreover, for $p\ge p_0$ and $c_0n\ge p$
 \[  \kla \intt h_1^p dP_{\mu} \mer^{\frac1p} \gl
 c_0 \pl  \frac{p}{1+\ln p} \pl.\]
In particular, the order of growth is best possible.
\end{exam}

{\bf Proof:} Since all the coefficients of $\al_{ij}$ are either
$0$ or $1$, it is clear that we count the number of events that
$j_k=1$. This yields  the first equality. For the second, we
consider $2\le t\le \frac{n}{2}$ and  $j<t\le j+1$ and deduce
from the proof of Lemma \ref{Herz} and Stirling's formula  that
  \for
   P_{\mu}(h_1>t)& \gl&  P_{\mu}(h_1=j) \lel {n \choose
 j} n^{-n} (n-1)^{n-j} \\
 &\gl& e^{-2} (4\pi)^{-\frac12} j^{-\frac12}
 \frac{n^n}{j^j(n-j)^{n-j}} \pl n^{-n} (n-1)^{n-j}
 \\
 &\gl& e^{-2} (4\pi)^{-\frac12} j^{-(j+\frac12)}  \gl
  e^{-2} (4\pi)^{-\frac12} (t-1)^{-(t-\frac12)}
 \pl.
 \mel
Therefore, we deduce from the proof of Lemma \ref{cc1} that for
$p\ge p_0$ and $\frac{n}{2}-1\ge \frac{4(p-\frac32)}{1+\ln
(p-\frac32)}$, we have
 \for
 \intt h_1^p dP_{\mu} &=& p \intt_0^\infty t^{p-1}
 P_{\mu}(h_1>t) dt \gl
 \frac{p}{e^2 (4\pi)^{\frac12}} \pl  \intt_2^{\frac n2}
 t^{p-1} (t-1)^{-(t-\frac12)} dt \\
 &=& \frac{p}{e^2 (4\pi)^{\frac12}} \pl  \intt_1^{\frac n2-1}
 (t+1)^{p-1} t^{-(t+\frac12)} dt \gl
  \frac{p}{e^2 (4\pi)^{\frac12}} \pl  \intt_1^{\frac n2-1}
 t^{p-\frac32} t^{-t} dt \\
 &\ge& \frac{p}{e^2 (4\pi)^{\frac12}} \pl e^{-2p-2}
 \kla \frac{p-\frac32}{1+\ln(p-\frac32)}
 \mer^{p-\frac12} \pl .
 \mel
This yields the assertion.\qed

\begin{rem} By complex interpolation, we see that for $1\le
q\le \infty$, $1\le p\le \infty$, we have
 \[ \kla n^{-n} \noo \summ_{j_1,...,j_n=1}^n
  \summ_{k=1}^n \al_{kj_k} e_k\rrm_{q}^p
  \mer^{\frac1p} \kl c_p^{\frac1q} \pl
 \noo  \summ_{k=1}^n \al^*_{(k-1)n+1}e_k \rrm_q \pl
 .\]
Again the same example shows that for $c_0p\ge q$ this behaviour
is best possible.
\end{rem}

\section{Application to noncommutative \boldmath$L_p$ \unboldmath-spaces}

This part is devoted to application in terms of non-commutative
version of symmetric spaces. Indeed, if $X$ is a symmetric
sequence space with basis $(e_k)$, one may define
 \[ S_X \lel \left\{ a\in B(\ell_2) \mitt  \noo \summ_k
 s_k(a)e_k\rrm_X\pl<\pl \infty \right\} \pl . \]
Here $(s_k(a))$ denotes the sequence of singular values, i.e. the
non-increasing rearrangement of the sequence of eigenvalues of
$\la_k((a^*a)^{\frac12})$. Then the norm of $a\in S_X$ is given by
 \[ \noo x\rrm_{S_X} \lel \noo \summ_k
 s_k(a)e_k\rrm_X \pl .\]
We refer to \cite[Proposition III.G.11]{W} and \cite[Proposition
2.a.5]{LTII}   to the non-trivial fact that this provides indeed
a norm, see \cite{DD,DDP} for more information. We use the
notation $S_X^n$ for the subspace of $n\times n$ matrices in
$S_X$. $M_n$ denotes the space of $n\times n$ matrices with the
operator norm. Let us also recall the more general definition of
noncommutative $L_p$ spaces. If $N$ is a von Neumann algebra and
$\tau$ is a normal, faithful, semifinite  trace, then the
$L_p$-norm of a $\tau$-measurable operator $x$ is defined as
 \[ \noo x\rrm_{L_p(N,\tau)} \lel [\tau(
 (x^*x)^{\frac p2})]^{\frac1p} \pl .\]
We refer to \cite{Ne} for basic properties and to \cite{KF} for
more information. The definition of $L_p$ spaces was extended to
non-semifinite von Neumann algebras by Connes \cite{Co} and
Haagerup \cite{Haa}. We only need the very basic properties
$b),c)$ (see \cite{Te,Pvp,JR})   and the recent result of Raynaud
that the class of non-commutative $L_p$ spaces is closed by
ultra-products (noted as a) below).

\begin{enumerate}
\item[a)] The class of non-commutative $L_p$ spaces
is closed by ultra-products.
\item[b)] $L_p(N)$ decomposes into $L_p(N)_{sa}+iL_p(N)_{sa}$ such that
 \[ \noo a+ib\rrm_p\sim_{2} \max\{\noo a\rrm_p,\noo
 b\rrm_p\}\pl .\]
\item[c)] For all $m\in \nz$, there is a distinctive norm $\noo \pl
\rrm_p$  on $L_p(M_m\ten N)$ such that for every unitary $u\in
M_m\bar{\ten}N$ and $x\in L_p(M_m\ten N)$
 \[ \noo u^*xu\rrm_{p} \lel \noo x\rrm_p \pl. \]
Moreover, for a diagonal n element $x=(x_{ij})\in L_p(M_m\ten N)$
(i.e. $x_{ij}=0$ for $i\neq j$), we have
  \[ \noo x\rrm_{L_p(M_m\ten N))} \lel \kla \summ_{i=1}^n \noo x_{ii} \rrm_{L_p(N)}^p \mer^{\frac1p}
  \pl .\]
\end{enumerate}

\begin{lemma} \label{sx} Let $1\le p< \infty$, $N$ be
a von Neumann algebra,   $X$ be a symmetric sequence space and
$\iota:X\to L_p(N)$ be an embedding. Then there is an embedding
of $S_X^n$ into $L_p(M_{n^n}\ten (N\oplus N))$.
\end{lemma}

{\bf Proof:} Let $v:X\to L_p(N)$ be an isomorphism into its
image, then we can decompose $v=(v_1,v_2)$ into two real-linear
maps $v_1(x)=\frac{v(x)+v(x)^*}{2}$ and
$v_2(x)=\frac{v(x)-v(x)^*}{2i}$ such that $\noo v(x)\rrm_p\sim_2
(\noo v_1(x)\rrm_p^p+\noo v_2(x)\rrm_p^p)^{\frac1p}$. Note that
the new map $(v_1,v_2):X\to L_p(N\oplus N)$ maps $X$ into the
selfadjoint part and is still a real-linear isomorphism. Thus, we
may assume that $v(X)\subset L_p(N)_{sa}$. Let $a\in M_n$ be a
selfadjoint matrix. Let $D_{\si}$ be the diagonal matrix given by
the sequence $\si=(s_k(a))_{k=1}^n$ and $u$ be unitary such that
$a=u^*D_{\si}u$. We deduce from Proposition \ref{matri} applied
to the matrix $\al_{ij}=\si_j(a)$.
 \for
 \noo a\rrm_X &=& \noo \summ_{k=1}^n
 s_k(a)e_k\rrm_X
  \pl \sim_{c_p} \pl  \kla n^{-n} \summ_{j_1,...,j_n=1}^n\noo \summ_{k=1}^n
 s_{j_k}e_k \rrm_X^p \mer^{\frac1p} \\
 &\sim_{c(v)}&  n^{-\frac{n}{p}} \kla\summ_{j_1,...,j_n=1}^n  \noo \summ_{k=1}^n
 s_{j_k}v(e_k) \rrm_{L_p(N)}^p \mer^{\frac1p}  \pl .
 \mel
Here $c(v)=\noo v\rrm \noo v^{-1}:Im(v)\to X\rrm$ only depends on
$v$.  For fixed $k$, we consider the map $\pi_k:M_n\ten M_{n^n}$
given by
 \[ \pi_k(x) \lel 1 \ten \cdots \ten
 \underbrace{x}_{\mbox{\scriptsize $k$-th
 position}} 1\ten \cdots \ten 1 \pl .\]
Note that $\pi_k(D_{\si})$ is a diagonal matrix $D_{\si_k}$ in
$M_{n^n}$ with entry $\si_k(j_1,...,j_n)=\si_{j_k}$. Therefore,
we have shown that for all diagonal matrices $x=D_{\si}$ we have
 \for
 \noo D_{\si} \rrm_X &\sim_{c_pc(v)} & n^{-\frac{n}{p}}
 \noo \summ_{k=1}^n \pi_k(D_{\si_k})\ten v(e_k)
 \rrm_{L_p(M_{n^n}\ten N)} \pl .
 \mel
However, $a=u^*D_{\si}u$ and
 \[ (u^*\ten \cdots \ten u^*)\pi_k(D_{\si})(u\ten \cdots
 \ten u)\lel \pi_k(a) \pl .\]
Since $(u\ten \cdots \ten u)\ten 1$ is a unitary in $M_{n^n}\ten
N$, we deduce
 \for
 \noo a \rrm_X &\sim_{c_pc(v)} & n^{-\frac{n}{p}}
 \noo \summ_{k=1}^n \pi_k(a)\ten v(e_k)
 \rrm_{L_p(M_{n^n}\ten N)} \pl .
 \mel
for all selfadjoint matrices $a \in M_n$ . Note that (since
$v(X)\subset L_p(N)_{sa})$ the map $T:S_X^n\to L_p(M_{n^n}\ten
N)$ defined by
  \[ T(a)\lel n^{-\frac{n}{p}} \pl \summ_{k=1}^n \pi_k(a)\ten v(e_k) \]
maps selfadjoint elements to selfadjoint elements. Thus for
arbitrary $x=a+ib$ we have
 \[ T(x)\lel T(a)+iT(b) \pl .\]
By $b)$ and the fact that
 \[ \noo x\rrm_{S_X} \sim_2 \max\{\noo a\rrm_{S_X},\noo
 b\rrm_{S_X}\} \]
we see that $T$ still defines an isomorphism with constant $16
c_p\noo v\rrm \noo v^{-1}:Im(v)\to X\rrm$.\qed

\begin{rem} The same remark shows that $S_q^n=S_{\ell_q}^n$ embeds
into the vector-valued non-commutative $L_p$-space
$S_p^{n^n}[\ell_q^n]$ defined by Pisier. Indeed,
 \[ \noo a\rrm_{S_q^n} \sim_{c_p} n^{-\frac{n}{p}} \pl \noo \summ_{k=1}^n \pi_k(a)\ten
 e_k \rrm_{S_p^{n^n}[\ell_q^n]} \pl .\]
We refer to \cite{Pvp} for a definition of the norm in this space
which is $\ell_p^{n^n}(\ell_q^n)$ on the diagonal and satisfies
 \[ \noo \summ_{k=1}^n ux_kv\ten e_k \rrm_{S_p^{n^n}[\ell_q^n]}
 \lel \noo \summ_{k=1}^n x_k \ten e_k \rrm_{S_p^{n^n}[\ell_q^n]}
 \]
for all unitaries $u,v\in M_{n^n}$. The embedding obtained in
this way is not a complete isomorphism for $q=1$.
\end{rem}

{\bf Proof:} Since $\noo a\rrm_1 \sim \max\{ \noo
\frac{a+a^*}{2}\rrm_q,\noo \frac{a-a^*}{2i}\rrm_q\}$, it suffices
to prove the equivalence of norms for selfadjoint matrices. Using
the unitary invariance, it suffices (as above) to prove it for
diagonal matrices. In that case it is a special case of
Proposition \ref{matri}. Now, let us indicate why this is not a
cb-isomorphism for $q=1$. We will freely use results from
\cite{Pvp}. Let us use  the notation $\tau_m=\frac{1}{m}tr$ for
the normalized trace. We consider the element
 \[ x\lel \summ_{i,j=1}^n e_{ij}\ten e_{ij} \in L_p(M_n,\tau_n;S_1^n) \pl
 .\]
Using simple facts about the Haagerup tensor product, we have
 \[ \noo x\rrm_{L_p(M_n,\tau_n;S_1^n)} \lel n^{-\frac1p} \noo id \rrm_{S_{p}^n[S_1^n]}^2 \lel 1
 \pl .\]
Due to the decomposition
 \begin{eqnarray}
 x &=&  (\summ_{i=1}^n e_{(i,i),1}) (\summ_{j=1}^n e_{1,(j,j)})
  \end{eqnarray}
we see that $x$ is positive. Here $e_{(ij),(k,l)}$ corresponds to
the matrix units in $M_{n^2}$. Positivity (see e.g. \cite{JD})
and the Burkholder/Rosenthal inequality (see \cite{JX}) imply
 \for  \lefteqn{
 \noo \summ_{k=1}^n id \ten \pi_k(x)\ten e_k
 \rrm_{L_p(M_{n^{n+1}},\tau_{n^{n+1}};\ell_1^n)}^{\frac12} \lel
  \noo \summ_{k=1}^n id \ten
 \pi_k(x)\rrm_{L_p(M_{n^{n+1}},\tau_{n^{n+1}})}^{\frac12} }\\
 & & \lel  \noo \summ_{k=1}^n (id \ten
 \pi_k(x^{\frac12})) \ten e_{k,1} \rrm_{L_{2p}(M_{n^{n+1}},\tau_{n^{n+1}};\ell_2^{c})}  \\
 & & \sim_{c_{2p}}  \max\{n^{\frac{1}{2p}} \noo
 x^{\frac12}\rrm_{L_{2p}(M_{n^2},\tau_{n^2})}, \sqrt{n} \noo
 \sum\nolimits_{i,j=1}^n e_{ij}\tau_n(e_{ij})
 \rrm_{L_p(M_n,\tau_n)}^{\frac12} \}\\
 & & \lel
 \max\{n^{\frac{1}{2p}} \noo
 x \rrm_{L_{p}(M_{n^2},\tau_{n^2})}^{\frac12} ,  \noo \sum\nolimits_{i=1}^n e_{ii} \rrm_{L_p(M_n,\tau_n)}^{\frac12}
 \} \\
  & & \lel \max\{ \noo x\rrm_{L_p(M_n,\tau_n;S_p^n)}^{\frac12},1\}  \pl .
 \mel
Now, we use the fact that the inclusion $id:S_1^n\to S_p^n$ is
not completely bounded and this is witnessed by the element  $x$.
Indeed, according to $(5)$, we see that $x$ represents a rank $1$
matrix and thus get
 \[ \noo  x \rrm_{L_{p}(M_{n},\tau_{n};S_p^n)} \lel n^{-\frac{1}{p}}
 \noo x\rrm_{S_p^{n^2}} \lel  n^{-\frac{1}{p}}
 \noo \summ_{i=1}^{n} e_{(i,i),1} \rrm_{2p}^2
 \lel n^{1-\frac{1}{p}}  \pl .\]
Therefore, the cb-norm of the map $T:S_1^n\to
L_p(M_{n^n},\tau_{n^n};\ell_1^n)$ satisfies
 \[ n^{1-\frac1p} \kl c_p \noo id_{L_p(M_n,\tau_n)}\ten T: L_p(M_n,\tau_n;S_1^n)\to
 L_p(M_{n^{n+1}},\tau_{n^{n+1}};\ell_1^n)\rrm \kl c_p \noo T\rrm_{cb}  \pl .\]
(Actually it is not very difficult to show that the upper
estimate holds too.) The assertion is proved.\qed

For the proof of Theorem \ref{m4}, we need some facts about
symmetric spaces with finite cotype.

\begin{lemma}\label{koe} Let $X$ be a symmetric sequence space
such that $X$ admits an embedding into $L_p(N)$ for some $1\le
p<\infty$, then for all $x\in X$
 \begin{eqnarray}
 \noo x\rrm_X &=& \sup_n \noo \summ_{k=1}^n x_ke_k\rrm_X  \pl .
 \end{eqnarray}
Moreover, the sequences with finite  support are dense.
\end{lemma}

{\bf Proof:} Since $L_p(N)$ has cotype $\max(2,p)<\infty$ (see
\cite{TJ,Fa}), we see that $X$ cannot contain a copy of $c_o$ on
disjoint blocks. Therefore $X$ is $\si$-order complete (cf
\cite[Proposition 1.a.5]{LTII}). From \cite[Proposition
1.a.7]{LTII}, we deduce that $X$ is $\si$-order continuous. Let
$x\ge 0$ be an element in $X$ and  consider
$y_n=x-\summ_{k=1}^nx_ke_k$. Then $y_n$ converges to $0$
everywhere and thus
 \[ \lim_n \noo y_n\rrm_X \lel 0 \pl .\]
This implies both assertions. \qed

\begin{lemma}\label{max}
Let $X\subset c_o$ be a symmetric sequence space satisfying the
Fatou property  $(6)$. For $k\in \nz$, we denote by $p_k$ the
projection onto the first $k$ unit vectors in $\ell_2$. Then
 \[ \noo a \rrm_{S_X} \lel \sup_k \noo p_kap_k
 \rrm_{X} \pl .\]
\end{lemma}

{\bf Proof:} Since $s_j(p_kap_k)\le s_j(a)$ we only have to show
$"\le"$. Since $X\subset c_o$, we may assume that $a$ is compact
(and using the spectral theorem for the compact operator
$(a^*a)^{\frac12}$) and thus, we may write $a=uD_{\si}v$ for
unitaries $u,v$ and a diagonal operator $D_{\si}$. In particular,
we may find projections $e_n\lel up_nu^*$ and $f_n=v^*p_nv$ of
rank $n$ such that
 \[ \noo up_nD_{\si}p_nv \rrm_{S_X} \lel \noo e_naf_n
 \rrm_{S_X} \pl .\]
Since $f_n$ and  $e_n$ have  finite ranks, we see that
 \for
 \lim_k \noo (1-p_k)f_n\rrm &=& 0  \lel \lim_k \noo e_n(1-p_k)\rrm \pl .
 \mel
By the triangle inequality, we deduce
 \for
 \lefteqn{  \noo e_naf_n\rrm_{S_X} }\\
  & & \lel  \lim_k \noo
  e_np_kap_kf_n+ e_n(1-p_k)ap_kf_n+
  e_np_ka(1-p_k)f_n + e_n(1-p_k)a(1-p_k)f_n\rrm \\
  & & \kl  \limsup_k \noo e_np_kap_kf_n\rrm_{S_X}+ \limsup_k 2 \noo
  e_n(1-p_k)\rrm_\infty \noo a\rrm_{S_X} + \noo
  a\rrm_{S_X} \noo (1-p_k)f_n\rrm_\infty \\
  & & \kl  \limsup_k \noo p_kap_k\rrm_{S_X} \pl .
  \mel
Since $X$ is supposed to satisfy $(6)$, we have
 \[\noo a\rrm \lel \sup_n \noo p_nD_{\si}p_n\rrm_{X}
 \kl \sup_n \noo e_naf_n\rrm_{S_X} \kl \sup_k \noo
 p_kap_k\rrm_{S_X} \pl . \]
The assertion is proved.\qed

{\bf Proof of Theorem \ref{m4}:} Since the space of diagonal
matrices in $S_X$ is $X$, it suffices to show that $S_X$ embeds
into some $L_p(N)$ if there is an isomorphism $v:X\to L_p(N)$.
Since $L_p(N)$ has cotype $\max(2,p)$, we have in particular, that
$X\subset c_o$. According to Lemma \ref{koe} and Lemma \ref{sx},
we see that
 \[ \noo a\rrm_{S_X} \lel \sup_n \noo p_nap_n\rrm_{S_X} \lel
 \limsup_n \noo p_nap_n\rrm_{S_X^n} \pl .\]
Let $\U$ be an ultrafilter on the integers. According to Lemma
\ref{max}, the mapping
 \[ \iota:S_{X}\to \prod_{n,\U} S_{X}^n  \quad,\quad \iota(x) \lel ((x_{ij})_{i,j=1}^n)  \]
is an isometric isomorphism. Due to Lemma \ref{sx}, we may find
$N_n$ and $T_n:S_X^n\to L_p(N_n)$ such that $\noo T_n\rrm\le 1$
and $\noo T_n:Im(T_n)\to S_X^{n}\rrm \le 16 c_p c(v)$. Hence,
 \[ T:S_X\to \prod_{n,\U} L_p(N_n) \]
is an isomorphism and the assertion is proved using Raynaud's
\cite{Ra} result (stated as a) above). \qed

\newcommand{\ww}{\vspace{-0.3cm}}

\end{document}